# The Investigation of Stationary Points in Central Configuration Dynamics

A.E. Rosaev [1]

[1] FGUP NPC "NEDRA" Yaroslavl, Russia

ABSTRACTS. A system of *N* points, each having *m* mass, and a central mass *M*, forming a planar central configuration, is considered. The motion equations for a testing particle are given in different co-ordinates. Different forms of motion equations are compared. For a higher *N*, the equations correct both for odd, and even number of particles are given. The linearization of equations is applied. The proposed method to calculate relevant sums can be important for different application.

It is shown, that the stationary solution (libration points) in considered system may be determined from algebraic equation of 5-th degree. It is obtained, that at large mass of particles outside libration point disappear. For inner libration points the limit minimal distance 1 was detected. In case small m / M solution for libration points in the considered problem have as a limit similar solution for 3-body collinear points.

The next step of our investigation is to construct another central configuration of 2N+1 body, where *N* new particles, each having $m_1$ mass, different from *m*, are placed to the libration point of our first system.

The main analytic results of this work verified with computer algebra system Maple® or by numeric calculations.

## 1. MAIN EQUATIONS

The system of N+1 particles is constructed, so one body placed at the center of the system, and another N points formed the regular polygon, vertexes of which placed on the described circle.

The stationary solution of N-body problem, at which masses placed in vertexes of the regular polygon is well known (Wintner,1941). The system of K similar central configurations is researched (Seidov,1990), where the equilibrium condition was given. The ability of capture in such system is shown (Pustylnikov, 1991), where the final motion considered. So, this model isn't new. Maxwell was first, who apply this model to the planetary rings dynamics. There are some evidences of presence large bodies in Saturn's ring obtained recently (Showalter, 1998). On the other side, as it known at present, the planetary rings is a very complex collisional system, which more suitable described by kinetic equations (Burns,1999).

A system of *N* points, each having $m_j$ mass, and a central mass *M*, forming a central configuration, is considered. The test particle is a particle of central configuration, motion of which is perturbed by another N-1 central configuration's particles. But it is easy to consider the case, when the test particle placed in the equal distance between two central configuration particles (in non-collinear libration points). The motion equations for a testing particle may be given in different co-ordinates. All of them can be derived by differentiation of potential energy. In accordance with (Ollongren,1981), the perturbation function in this system is:

$$U = \frac{GM}{R} + \sum_{j=0}^{N-1} Gm_j \frac{1}{\left((X-X_j)^2 + (Y-Y_j)^2 + (Z-Z_j)^2\right)^{1/2}} \qquad (1.1)$$

where G -gravity constant, X,Y,Z, $X_j$, $Y_j$, $Z_j$ - rectangular coordinates of testing particle m and particles of central configuration $m_j$ accordingly, R- distance between M and m.

The equations of motion of central configuration's particle may be obtained by the calculation of the partial derivatives of U, the equations for n central configurations easy followed from one central configuration case.

Ollongren (1981) and some other authors use different expressions of perturbation function for odd and even number of central configuration's particles in rectangular coordinates. Perov in (Rosaev, Perov, Sadovnikov, 1996) suggest very complicated expressions, using oblatness of central body in additions to mutual central configuration's particles perturbations.

In our opinion, the motion equations in the polar coordinates are more suitable to study central configurations. They can be derived by common way – the differentiation (1.1). Now, we can consider their derivation immediately from the task's geometry.

We shall consider planar motion of particle with mass m in field of attraction other N-1 particles, placed in vertexes of regular polygons and central body with mass M. Evidently, for radial and tangential forces of interaction between particles we have:

$$F_r = -\sum_j G m m_j \frac{\cos\varphi_j}{r_j^2} \qquad F_{tg} = -\sum_j G m m_j \frac{\sin\varphi_j}{r_j^2} \qquad (1.2)$$

where $r_j$ - distance between particles, $\varphi_j$ - angle between direction $Om$ and $mm_j$, $O$ – the configuration's center. In case of circular central configuration:

$$r_j^2 = x^2 + (2R\sin(\alpha_j/2))^2 (1 + x/R)$$
$$\cos\varphi_j = \frac{2Rx + 2x^2 + (2R\sin(\alpha_j/2))^2 (1 + x/R)}{2(R+x)\sqrt{x^2 + (2R\sin(\alpha_j/2))^2 (1 + x/R)}} \qquad (1.3)$$

where $R$ – the central configuration's radius, $\alpha_j$ - angle between particles, and $x$- a distance between the test particle and the central configuration. The motion equation can easily be obtained by replacement of the expressions (1.3) to (1.2) ones. One can see, that these equations are correct both for odd, and even number of particles:

$$\ddot{R} - R(\dot\varphi + \Omega)^2 = -\frac{GM}{R^2} - \sum_{j=1}^{N-1} Gm_j \frac{2R_0 \sin^2(\alpha_j/2) + x}{\left(x^2 + (2R_0\sin(\alpha_j/2))^2(1 + x/R_0)\right)^{3/2}}$$

$$R\ddot\varphi + 2(\dot\varphi + \Omega)\dot x = -\sum_{j=1}^{N-1} Gm_j \frac{2R_0 \sin(\alpha_j/2)\cos(\alpha_j/2)}{\left(x^2 + (2R_0\sin(\alpha_j/2))^2(1 + x/R_0)\right)^{3/2}}$$

(1.4)

where $R = R_0 + x$ - the distance of testing particle from the center, $R_0$ - the distance of $j$-particle from the center, $\Omega$ - angular velocity, $\alpha_j = 2\pi\, j/N + \varphi$, where $\varphi$ is possible angular declination from stationary position.

Now, we derive the motion equation by common way. In cylindrical system:

$$m_i \left( \ddot R_i - R_i \dot\lambda_i^2 \right) = \frac{dU}{dR_i} \equiv U_r$$

$$m_i \frac{d(R_i^2 \dot{\lambda}_i)}{dt} = \frac{dU}{d\lambda_i} \equiv U_l \quad (1.5)$$

$$m_i \frac{d^2 z_i}{dt^2} = \frac{dU}{dz} \equiv U_z$$

The perturbation function in cylindrical frame:

$$U := \sum_{i=1}^{N-1} \frac{G m m_i}{\sqrt{r(i)^2 + R^2 - 2 r(i) R \cos(\alpha) + (z - z(i))^2}}$$

Here, for simplification, $\alpha = \alpha_i$. After differentiation and substitution $R = R_0 + x(t)$; $r(i) = R_0 \equiv R$:

$$U_l := \sum_{i=1}^{N-1} \left( -\frac{G m m_i (R+x) R \sin(\alpha)}{\left((R+x)^2 + R^2 - 2(R+x) R \cos(\alpha) + (z - z(i))^2\right)^{3/2}} \right)$$

$$U_r := \sum_{i=1}^{N-1} \left( -\frac{1}{2} \frac{G m m_i (2R + 2x - 2R \cos(\alpha))}{\left((R+x)^2 + R^2 - 2(R+x) R \cos(\alpha) + (z - z(i))^2\right)^{3/2}} \right)$$

(1.6)

These equations were derived with the use of Maple® software for symbolic computations. After insignificant manual simplifications, the motion equations can be rewritten in the following way (1.4). Evidently, for a great *N*, these equations are correct both for odd, and even number of particles.

Let us assume the masses of all particles to be equal, $m_j = m$. In this case the second equation in (1.4) reduces to conservation of the angular moment:

$$L = (\dot{\varphi} + \Omega) R^2 = const \quad (1.7)$$

Then, to solve this problem, we try to use some existing software tools. Only Maple® allows to get a solution in symbolic form, and to consider different asymptotic cases. Different solutions for this equation are compared. As before, we assume *x* to be the deviation of the particle from the central configuration. The symbolic computations allow us to represent the solution as an *x* power series (in case of small *x*). This enables to study dependence of libration point position on initial conditions.

In case of insignificant perturbations, x<<R, the first equation (1.4) reduces to the form (*L* - angular moment):

$$\ddot{x} - \frac{2Gm}{R^3} x + \frac{3L^2}{R^4} x = -\sum_j \frac{Gm}{(2R \sin(\alpha_j / 2))^3} x \quad (1.8)$$

The expansion (1.6) by power x, obtained with Maple, proves this result:

$$U_r := \left( \sum_{i=1}^{N-1} \left( -\frac{1}{2} \frac{G m m_i (2R - 2R\cos(\alpha))}{\left(2R^2 + (z-z(i))^2 - 2R^2\cos(\alpha)\right)^{3/2}} \right) \right) + \left( \sum_{i=1}^{N-1} \left( -\frac{G m m_i}{\left(2R^2 + (z-z(i))^2 - 2R^2\cos(\alpha)\right)^{3/2}} + \frac{3}{4} \frac{G m m_i (2R - 2R\cos(\alpha))^2}{\left(2R^2 + (z-z(i))^2 - 2R^2\cos(\alpha)\right)^{5/2}} \right) \right) x + \left( \sum_{i=1}^{N-1} \left( \frac{3}{2} \frac{G m m_i (2R - 2R\cos(\alpha))}{\left(2R^2 + (z-z(i))^2 - 2R^2\cos(\alpha)\right)^{5/2}} - \frac{1}{2} G m m_i \right. \right.$$

$$\left. \left. \left( -\frac{3}{2} \frac{1}{2R^2 + (z-z(i))^2 - 2R^2\cos(\alpha)} + \frac{15}{8} \frac{(2R - 2R\cos(\alpha))^2}{\left(2R^2 + (z-z(i))^2 - 2R^2\cos(\alpha)\right)^2} \right) \right.\right.$$

$$\left.\left. (2R - 2R\cos(\alpha)) \Big/ \left(2R^2 + (z-z(i))^2 - 2R^2\cos(\alpha)\right)^{3/2} \right) \right) x^2 + O(x^3)$$

(1.9)

$$U_l := \left( \sum_{i=1}^{N-1} \left( -\frac{G m m_i R^2 \sin(\alpha)}{\left(2R^2 + (z-z(i))^2 - 2R^2\cos(\alpha)\right)^{3/2}} \right) \right) + \left( \sum_{i=1}^{N-1} \left( -\left( \frac{G m m_i}{\left(2R^2 + (z-z(i))^2 - 2R^2\cos(\alpha)\right)^{3/2}} - \frac{3}{2} \frac{G m m_i (2R - 2R\cos(\alpha)) R}{\left(2R^2 + (z-z(i))^2 - 2R^2\cos(\alpha)\right)^{5/2}} \right) R \sin(\alpha) \right) \right) x + \left( \sum_{i=1}^{N-1} \left( \left( -\frac{3}{2} \frac{G m m_i (2R - 2R\cos(\alpha))}{\left(2R^2 + (z-z(i))^2 - 2R^2\cos(\alpha)\right)^{5/2}} + G m m_i \right. \right. \right.$$

$$\left. \left. \left. \left( -\frac{3}{2} \frac{1}{2R^2 + (z-z(i))^2 - 2R^2\cos(\alpha)} + \frac{15}{8} \frac{(2R - 2R\cos(\alpha))^2}{\left(2R^2 + (z-z(i))^2 - 2R^2\cos(\alpha)\right)^2} \right) R \right. \right. \right.$$

$$\left. \left. \left. \Big/ \left(2R^2 + (z-z(i))^2 - 2R^2\cos(\alpha)\right)^{3/2} \right) R \sin(\alpha) \right) \right) x^2 + O(x^3)$$

(1.10)

It seems that close to equilibrium radial force depend from x, $U_r = a0+(a1+a2)x+a3*x^2+\ldots$ In the coefficient before x, the second term (a2) much smaller, then first (a1, table 1). Close to equilibrium position, the radial perturbation is much more than tangential (Fig.1).

Table 1

The structure of coefficient before x in the radial perturbation expansion

| N | a1 | a2 | a1/a2 |
|---|---|---|---|
| 10 | 39.88 | 7.26 | 5.5 |
| 100 | 38789 | 129.3 | 300 |
| 1000 | $3.877*10^7$ | 1844.4 | $2.1019*10^4$ |
| 2000 | $3.102*10^8$ | 4020 | $7.7164*10^4$ |
| 4000 | $24.804*10^8$ | 8701 | $28.5068*10^4$ |
| 6000 | $83.790*10^8$ | 13635 | $61.4520*10^4$ |
| 8000 | $198.80*10^8$ | 18733 | $106.1340*10^4$ |
| 10000 | $386.17*10^8$ | 23919 | $161.4500*10^4$ |

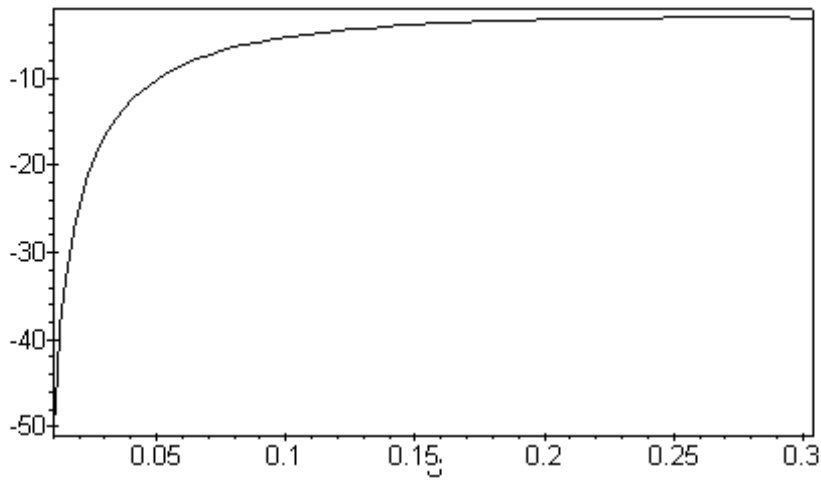

Fig.1 Radial/tangential force ratio (in dependence from S)

As a result, after the substitution (1.9) in (1.5) we have equation (1.8). So, in linear approximation, it is oscillations with a frequency of:

$$\omega = \sqrt{-\frac{2Gm}{R^3} + \frac{3L^2}{R^4} + \sum_j \frac{Gm}{(2R\sin(\alpha_j/2))^3}} \quad (1.11)$$

The calculation $U_l$ according (1.10), shows, that tangential force much smaller then radial. The second term of expansion (1.9) gives the conditions appearance of oscillations:

$$\sum_j m_j/(2R\sin(\alpha_j/2))^3 \gg 6M/R^4 x \quad (1.12)$$

*Rosaev A.E. The investigation of stationary points in central configuration dynamics*

The perturbations, caused by moonlet with mass *m,* move the ring's particles: near ones first, and then distant ones. As the ring is not absolutely rigid, the phase angle δ will appear, and in general case:

$$\ddot{x} = -\omega^2 x - 2k\dot{x} + fx$$

$$tg\delta = \frac{2k\omega}{\omega_0^2 - \omega^2} \tag{1.13}$$

with $k$ – resistance factor, $\omega_0 = \sqrt{-\frac{2Gm}{R^3} + \frac{3L^2}{R^4}}$ - fundamental frequency. The difference between ω and $\omega_0$ enables to describe stationary oscillations (Fig.2) within the central configuration. Their existence is conditioned by:

$$\frac{2\pi}{\lambda} = \frac{\omega}{\omega_0} = \sqrt{\frac{\sum_j G m_j /(2R \sin(\alpha_j/2))^3}{-\frac{2Gm}{R^3} + \frac{3L^2}{R^4}}} + 1 = \gamma \tag{1.14}$$

where γ - integer. The character of motion is conditioned by *k* value. In case of high *k* the oscillations damp quickly, forming a non-vibrating component.

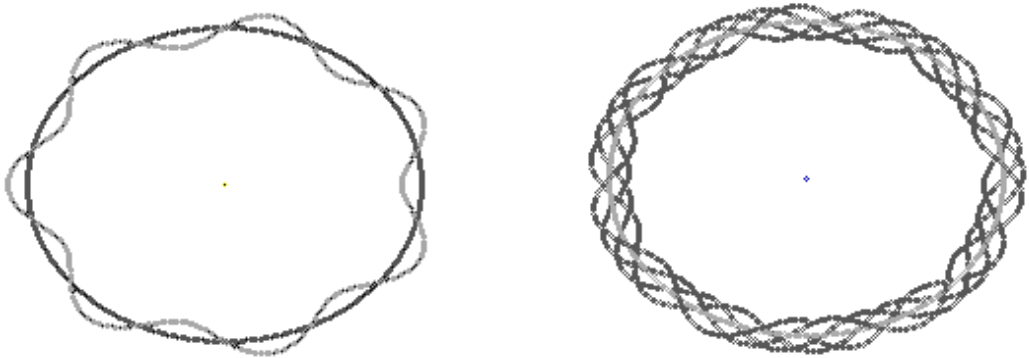

Fig.2. Stationary and non-stationary oscillations

On the whole, the time of existence of this kind of motion is determined, by collisions with non-vibrating component and may be long. One of possible reasons of these oscillations is a collision with the nearby moonlet.

From $R_0$=R=const followed, that dependent of the time part of angular velocity $\dot{\varphi}$ =0 and first equation (1.4) reduced to:

$$\ddot{R} - R_0\Omega^2 = -\frac{GM}{R_0^2} - \sum_{j=1}^{N-1} Gm \frac{1}{4R^2_0 \sin(\alpha_j/2)} \tag{1.15}$$

In case of the regular rotation of all system, $\ddot{R}$ = 0 and may to determine Ω. (for non-interacted particles evidently $\Omega_0^2$=GM/$R_0^3$):

$$\Omega = 1/\sqrt{R_0}\left(\frac{GM}{R_0^2} + \sum_{j=1}^{N-1} Gm \frac{1}{4R_0^2 \sin(\alpha_j/2)}\right)^{1/2} = \Omega_0\left(1 + \frac{m}{4M}\sum_{j=1}^{N-1}\sin^{-1}(\alpha_j/2)\right)^{1/2} \quad (1.16)$$

By using the approximate formula:

$$\sum_{j=1}^{N-1}\sin^{-1}(\alpha_j/2) = \frac{2N}{\pi}\left(\ln\frac{2N}{\pi} - 0.58\right) \quad (1.17)$$

may estimate the increasing of the angular velocity of central configuration at the different number of particles. By express the mass of central configuration in the central mass units, may be construct the dependence $\Omega(N)$. So, the accounting of the ring's particles attraction followed to moving resonance position away from planet. This displacement increases with the mass of ring and the number of particles increasing. At $N=10^{12}$ and $m=0.01\,M$ deposition reach 2%. For the small mass of the ring this effect is neglectable.

## 2. THE LIBRATION POINTS COORDINATES

To consider any questions, related with stability in this problem, first of all, the libration points should be determined. In this part of work we shell study only planar central configurations. The libration points coordinates determined from (1.4) at the requirements $\ddot{R} = 0, \ddot{\varphi} = 0, \dot{R} = 0, \dot{\varphi} = 0$. There are two kinds of solution exist:

$$\alpha_j = 2\pi\, j/N + \varphi_0$$

$$-R\Omega^2 = -\frac{GM}{R^2} - \sum_{j=1}^{N} Gm \frac{2R_0 \sin^2(\alpha_j/2) + x_0}{\left(x_0^2 + (2R_0 \sin(\alpha_j/2))^2(1 + x_0/R_0)\right)^{3/2}} \quad (2.1)$$

and

$$-R\Omega^2 = -\frac{GM}{R^2} - \frac{Gm}{x_0^2} - \sum_{j=1}^{N-1} Gm \frac{2R_0 \sin^2(\alpha_j/2) + x_0}{\left(x_0^2 + (2R_0 \sin(\alpha_j/2))^2(1 + x_0/R_0)\right)^{3/2}} \quad (2.2)$$

$$\alpha_j = 2\pi\, j/N + \varphi_0$$

In the first case $\varphi_0 = \pi/N$, $x_0=0$ and we have noncollinear libration point. In the second case $\varphi_0 = 0$, $x_0 \neq 0$ and x-coordinate of collinear points required to be determined. So, the equations (2.2) for the collinear points of libration are most interesting. It may be reduced to fifth degree polynom.

At small $x_0 \ll R_0$, approximately:

$$\sum_{j=1}^{N-1} Gm \frac{2R_0 \sin^2(\alpha_j/2) + x_0}{\left(x_0^2 + (2R_0 \sin(\alpha_j/2))^2(1 + x_0/R_0)\right)^{3/2}} \approx \sum_{j=1}^{N-1} Gm \frac{2R_0 \sin^2(\alpha_j/2) + x_0}{(2R_0 \sin(\alpha_j/2))^3} \quad (2.3)$$

After the substitution $R = R_0 + x_0$ and multiplying all terms to $x_0^2(R_0 + x_0)^2$:

$$-(R_0 + x_0)^3 x_0^2 \Omega^2 = -GMx_0^2 - Gm(R_0 + x_0)^2 - (Ax + B)(R_0 + x_0)^2 x_0^2 \qquad (2.4)$$

where:

$$A = \sum_{j=1}^{N-1} Gm \frac{1}{(2R_0 \sin(\alpha_j/2))^3}; \qquad B = \sum_{j=1}^{N-1} Gm \frac{1}{(4R_0^2 \sin(\alpha_j/2))} \qquad (2.5)$$

Put $x_0 = x$, $R_0 = R$. After the simplification:

$$(\Omega^2 + A)x^5 + (3\Omega^2 R + 2AR + B)x^4 + (3\Omega^2 R^2 + R^2 A + 2BR)x^3 + (BR \pm Gm)x^2 \pm GmR = 0 \qquad (2.6)$$

There are few ability to continue simplifications. The coefficient B is very small, due to sinus function symmetry, so coefficient proportional $x^2$ in (2.6) neglectable in all cases. If R large and x small,

$$x \ll 1 \ll R \qquad (2.8)$$

we can write approximately:

$$(3GM/R + AR^2 + 2BR)x^3 \approx \pm GmR^2 \qquad (2.8)$$

In general case, by taking into accounting, that $\Omega^2 \sim GM/R_0^3$, approximately solution of them (at $m/M \ll 7.84048/N^3$ and $x/R < \pi/N$) is:

$$x_L^3 \approx \pm m / (M/R^3(3-k)); \qquad k = 2.4041 \, mN^3 / (8M\pi^3) \qquad (2.9)$$

At the limit $M \gg m$, the expression (2.4) reduced to well-known solution 3-body problem (Szebehely,1967). At the increasing of masses of central configuration's particles, the moving of libration points away from central configuration's axis take place (table 1). The similar behavior observed, when the number of particles decreased. In neighborhood k=3 and at largest k expressions (2.4) are unsuitable.

Numerical solution at large m/M showed, that the inner libration point coordinate have a limit s, depended from N. It cannot close to the center of system near then s, even if M=0. Outer point of libration going away from the central configuration's axis for all time increasing m/M ratio, but at large mass of central configuration (m*N/M>1) outer libration points are disappeared (tables 2-3, where X0/R - 3-body solution). There are mass of central configuration (in central mass units) in the left column, and three pairs of columns with libration points coordinates (in according with number of particles N in the first row): $X_L/R$ - for N+1 bodies case, X0/R - for three bodies case.

As it seems from table 1, for inner points we have stop in changing of libration point co-ordinate $X_L/R$ for N+1 bodies since m*N/M>100, in contrast of the libration point co-ordinate

X0/R (three bodies), which continue to increase. For outer libration points, we have asymptotic increasing coordinates $X_L/R$ at range $0.1 < m*N/M < 1$ (table 3).

Table 2

Dependence of libration point's from mass ratio (inside point)

| m*N/ M | N=50 | | N=100 | | N=1000 | |
|---|---|---|---|---|---|---|
| | $X_L/R$ | X0/R | $X_L/R$ | X0/R | $X_L/R$ | X0/R |
| 0.00001 | -0.0045 | -0.0041 | -0.0035 | -0.0032 | -0.0015 | -0.0015 |
| 0.00010 | -0.0085 | -0.0087 | -0.0065 | -0.0069 | -0.0035 | -0.0032 |
| 0.00100 | -0.0185 | -0.0188 | -0.0145 | -0.0149 | -0.0095 | -0.0069 |
| 0.01000 | -0.0375 | -0.0405 | -0.0305 | -0.0322 | -0.0255 | -0.0149 |
| 0.10000 | -0.0705 | -0.0874 | -0.0605 | -0.0693 | -0.0505 | -0.0322 |
| 1.0 | -0.1055 | -0.1882 | -0.0915 | -0.1494 | -0.0655 | -0.0693 |
| 10.0 | -0.1175 | -0.4055 | -0.0995 | -0.3218 | -0.0675 | -0.1494 |
| 100.0 | -0.1195 | -0.8736 | -0.1005 | -0.6934 | -0.0685 | -0.3218 |
| 1000.0 | -0.1195 | -1.8821 | -0.1015 | -1.4938 | -0.0685 | -0.6934 |
| 10000.0 | -0.1195 | -4.0548 | -0.1015 | -3.2183 | -0.0685 | -1.4938 |

Table 3

Dependence of libration point's from mass ratio (outside point)

| m*N/ M | N=50 | | N=100 | | N=1000 | |
|---|---|---|---|---|---|---|
| | $X_L/R$ | X0/R | $X_L/R$ | X0/R | $X_L/R$ | X0/R |
| 0.00001 | 0.0045 | 0.0041 | 0.0035 | 0.0032 | 0.0015 | 0.0015 |
| 0.00010 | 0.0085 | 0.0087 | 0.0075 | 0.0069 | 0.0035 | 0.0032 |
| 0.00100 | 0.0175 | 0.0188 | 0.0155 | 0.0149 | 0.0125 | 0.0069 |
| 0.01000 | 0.0475 | 0.0405 | 0.0425 | 0.0322 | 0.0425 | 0.0149 |
| 0.10000 | 0.1925 | 0.0874 | 0.2125 | 0.0693 | 0.2975 | 0.0322 |
| 1.0 | | | | | | |
| 10.0 | The libration points are absent | | | | | |
| 100.0 | | | | | | |

The result (2.9), obtained analytically, then proved with Maple® as a solution of (2.6) at limit B→ 0 (and small x):

$$x := \left(\frac{1}{R(A+3W)}\right)^{1/3} m^{1/3} + O\!\left(m^{2/3}\right), I\sqrt{R} + O\!\left(R^{3/2}\right)$$

(2.10)

Here W≡Ω, $I = \sqrt{-1}$. In addition, some another asymptotic expansions (2.6) was obtained with Maple. At limit m→ 0:

$$x3 := 0,\ -\frac{1}{R(A+3W)} B + O\!\left(B^2\right), O\!\left(R^{-1}\right)$$

(2.11)

### 3. GENERALIZATION FOR THE CASE 2N+1 PARTICLE

Let us consider the system of 2N masses $m_l$, l=1...n, which form two symmetric polygons, the vertex of each lying on the radius the circle described. There are two main cases (Fig 1 -2).

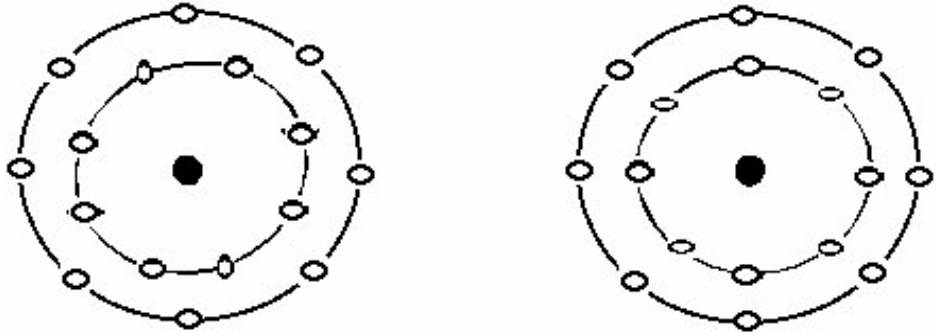

Fig.3 Noncollinear 2N+1 body configuration. (left) collinear 2N+1 body configuration (right)

The mass $M$ is in the system centre. Each pair of masses $m_l$ and $m_j$ lies on the same ray beginning in the system centre - collinear 2N+1 body configuration, or in a middle position - - noncollinear 2N+1 body configuration (Fig 4-5). The motion equations for particles of central configurations may be derived from equation (Rosaev, 1997) for one central configuration. In polar system we have:

$$x_k'' - (R_k + x_k)(\varphi_k' + \Omega_k)^2 = -\frac{GM}{(R_k + x_k)^2} + \widetilde{F}_k \pm \widetilde{F}_s$$

$$(R_k + x_k)\cdot\varphi_k'' + 2(\varphi_k' + \Omega_k)x_k' = f_i + f_o$$

$$k = o, i; \quad s = i, o$$

(3.1)

where $k = o, i$ ; $s = i, o$. There are four cases in these equations. But to study stability only two of them are of interest. And we need to choose "-" in the right part of the expressions to describe motion between rings.

$$\widetilde{F}_k = -\sum_{j=1}^{N-1} Gm_k \frac{2R_k \sin^2(\alpha_j/2) + x_k}{\left(x_k^2 + (2R_k \sin(\alpha_j/2))^2(1 + x_k/R_k)\right)^{3/2}}$$

$$f_k = -\sum_{j=1}^{N-1} Gm_k \frac{2R_k \sin(\alpha_j/2)\cos(\alpha_j/2)}{\left(x_k^2 + (2R_k \sin(\alpha_j/2))^2(1 + x_k/R_k)\right)^{3/2}}$$

(3.2)

$$\alpha_j = 2\pi\, j/N + \varphi$$

In these expressions $R_0$ is the central configuration radius, $N$ - number of particles, $m_k$ - mass of a central configuration's particle, $G$ - constant of gravity, $x_k$ and $\varphi$ - deviations of the testing particle from the equilibrium point in radial and tangential directions respectively, $\Omega$ - angular velocity of circular rotation. Evidently, the equation (3.1) is correct at any value of $N$, both odd and even, like expression (1.4). We assume, that $N$ is rather high.

Then, we need to use the expansion:

$$\tilde{F}_k = B_k + A_k x_k + O(x^2)$$

$$f_k = b_k + a_k x_k + O(x^2)$$

(3.3)

Here the coefficients should be obtained by $f$ and $F$ differentiation, and be represented by sums from $j=1$ to $N-1$:

$$A_k = -\frac{Gm_k}{R^3} \sum_{j=1}^{N-1} \left[ \frac{1}{|2\sin(\alpha_j/2)|^3} - \frac{3}{|8\sin(\alpha_j/2)|} \right]$$

$$B_k = -\frac{Gm_k}{4R^2} \sum_{j=1}^{N-1} \frac{1}{|\sin(2\pi\, j/N)|}$$

$$a_k = -\sum_{j=1}^{N-1} Gm_k \frac{-3(\sin(\alpha_j/2)\cos(\alpha_j/2))}{|(2R\sin(\alpha_j/2))^3|} = 0$$

$$b_k = -\sum_{j=1}^{N-1} Gm_k \frac{2R\sin(\alpha_j/2)\cos(\alpha_j/2)}{|(2R\sin(\alpha_j/2))^3|} = 0$$

$$\Omega_k^2 = \frac{GM}{R_k^3} + \sum_{j=1}^{N-1} \frac{Gm_k}{4R_k^3 |\sin(\alpha_j/2)|}$$

(3.4)

It seems, that coefficient $A$ is of the highest value among them, and it can allow numerical approximation:

$$A_k \approx -\alpha(N)\, Gm_k\, N^3/(2\pi R)^3\ ; \qquad \alpha(N) \approx 2.404101 \qquad (3.5)$$

the validity of which, as well as validity of linearization are shown in table 4. Besides, the expansion coefficients $\alpha'$ for $1/j^3$ sums are given.

Table 4

Coefficients of expansion depending on number of particles

| N | α(N) | α'(N) |
|---|---|---|
| 10 | 2.6580862 | 2.37132406 |
| 20 | 2.4808764 | 2.39506411 |
| 50 | 2.4196434 | 2.40257668 |
| 100 | 2.4086496 | 2.40372180 |
| 200 | 2.4054148 | 2.40401482 |
| 1250 | 2.4041442 | 2.40410137 |
| 2500 | 2.4041123 | 2.40410137 |
| 5000 | 2.4041037 | 2.40410137 |
| 10000 | 2.4041018 | 2.40410137 |
| 20000 | 2.4041016 | 2.40410137 |

As it seems from a Table 4,

$$\sum_{i=1}^{N} 1/(\sin(\pi i / N))^3 \longrightarrow N^3/(2\pi)^3 \sum_{i=1}^{N} 1/i^3 \longrightarrow 2.40410...N^3/(2\pi)^3$$

when N increase. It means, that only nearest neighbor particles give main contribution in perturbation.

This is a remarkable result, which makes it possible to effectively calculate perturbation effect in close neighborhood of the central configuration. This method was successfully applied above to determine coordinates of the libration points.

It is easy to see the limit of linear approximation: the error increases with *x/R* growth. The range of the linear approximation applicability decreases with the number of particles increase.

## 4. RESULTS AND CONCLUSIONS

Our final target is to construct the stable kN+1 points central configurations. To solve this problem, first of all, the stationary (libration) points in N+1 case need to be determined.

At the beginning, we considered a system of *N* points, each having *m* mass, and a central mass *M*, forming a central configuration. The motion equations for a general case (odd and even number of particles) are obtained by differentiation of potential function. The proposed method to calculate mutual perturbation force, obtained in this work, can be important for practical application.

Co-ordinates for the libration points in the system are calculated. It is shown, that the libration points in considered system may be determined from algebraic equation of 5-th degree. Main attention is given to the case of large N. It is obtained, that at large mass of particles outside libration point disappear. For inner libration points the limit distance l exist - libration points cannot be close then l to the central mass M at fixed N and arbitrary m / M ratio.

We used analytical approach as well as computer algebra method to study central configuration dynamics. As a result, we can conclude, that to describe libration points behavior at different number (*N*) of particles, analytical method is more preferable; in case of dependence on a particle mass we obtain the same results using the both methods. The result coincides with classical solution within the limit for the 3-body problem. Nevertheless, some results can be obtained most easily only with computer algebra system. In particular, asymptotic solution within the limit of central configuration's particle $m \to 0$ of a neglectable mass may be obtained only by

*Rosaev A.E. The investigation of stationary points in central configuration dynamics*

this way. The computer algebra methods have great advantage to expand the perturbation function in the described problem. By using Maple® tools, we can determine frequency of insignificant oscillations in the system.

In all these cases, computer algebra methods can help while solving complex problems, making the solutions easier, as well as while verifying the results. We believe, that Maple® computer algebra package is a suitable tool for this purpose. But it does not mean, that special computer packages to analyse dynamic systems are not required.